\documentclass[oneside,11pt]{amsart}
\usepackage{amsmath}
\usepackage{amssymb}
\usepackage{graphicx}
\usepackage{tikz}
\usepackage{color}
\usepackage[colorlinks=true]{hyperref}
\usepackage{tikz-cd}
\usepackage{todonotes}
\usepackage{thm-restate}
\usepackage{fullpage}
\newtheorem{thm}{Theorem}[section]
\newtheorem{prop}[thm]{Proposition}
\newtheorem{lem}[thm]{Lemma}

\theoremstyle{definition}
\newtheorem{defn}[thm]{Definition}
\newtheorem{rem}[thm]{Remark}

\newtheorem{exmp}[thm]{Example}
\newtheorem{ques}[thm]{Question}

\renewcommand{\bar}[1]{\overline{#1}}

\renewcommand{\emptyset}{\varnothing}

\newcommand{\field}[1]{\mathbb{#1}}
\newcommand{\Z}{\field{Z}}
\newcommand{\R}{\field{R}}

\newcommand{\PP}{\field{P}}

\newcommand{\cH}{\mathcal{H}}
\newcommand{\cG}{\mathcal{G}}

\newcommand{\cS}{\mathcal{S}}

\DeclareMathOperator{\lk}{Link}

\usepackage{ifthen}

\newcommand{\showcomments}{yes}

\newsavebox{\commentbox}
{\ifthenelse{\equal{\showcomments}{yes}}%
{\footnotemark
    \begin{lrbox}{\commentbox}
    \begin{minipage}[t]{1.25in}\raggedright\sffamily\upshape\tiny
    \footnotemark[\arabic{footnote}]}%
{\begin{lrbox}{\commentbox}}}%
{\ifthenelse{\equal{\showcomments}{yes}}%
{\end{minipage}\end{lrbox}\marginpar{\usebox{\commentbox}}}%
{\end{lrbox}}}

\title{Coarse geometry of extended admissible groups}

\author{Toan Trong Dao}
\address{Faculty of Mathematics and Computer Science, University of Science, Ho Chi Minh City, Vietnam}
\address{Vietnam National University,
Ho Chi Minh City, Vietnam}
\email{daotrongtoan.dtt4@gmail.com}

\author{Hoang Thanh Nguyen}
\address{FPT University, Da Nang, Vietnam}
\email{nthoang.math@gmail.com}

\date{\today}
 % Hyperbolic groups and nonpositively curved groups

\begin{document}

\begin{abstract}
Extended admissible groups belong to a particular class of graphs of groups that admit a decomposition generalizing those of non-geometric 3-manifold groups and Croke-Kleiner admissible groups. In this paper, we study several coarse-geometric aspects of extended admissible groups. We show that changing the gluing edge isomorphisms does not affect the quasi-isometry type of these groups. We also prove that, under mild conditions on the vertex groups, extended admissible groups exhibit large-scale nonpositive curvature, thereby answering a question posed in \cite{NY23}.

As an application, our results enlarge the class of extended admissible groups known to admit well-defined quasi-redirecting boundaries, a notion recently introduced by Qing–Rafi. In addition, we compute the divergence of extended admissible groups, generalizing a result of Gersten from non-geometric 3-manifold groups to this broader setting. Finally, we study several aspects of subgroup structure in extended admissible groups.
\end{abstract}

\maketitle

\section{Introduction}
Let $M$ be a non-geometric 3-manifold. The torus decomposition of $M$ yields a nonempty minimal union $\mathcal{T} \subset M$ of disjoint essential tori, unique up to isotopy, such that each component $M_v$ of $M \backslash \mathcal{T}$, called a \emph{piece}, is either  Seifert fibered or hyperbolic.
There is an induced graph of groups decomposition $\mathcal{G}$ of $\pi_1(M)$ with underlying graph $\Gamma$ as follows. For each piece $M_v$, there is a vertex $v$ of $\Gamma$  with vertex group $\pi_1(M_{v})$. For each torus $T_e\in \mathcal T$ contained in the closure of pieces $M_v$ and $M_{w}$, there is an edge $e$ of $\Gamma$ between vertices $v$ and $w$. The associated edge group is $\pi_1(T_e)\cong \Z^2$ and the edge monomorphisms are the maps induced by inclusion. 
Each Seifert fibered piece $M_v$ in the JSJ decomposition of $M$ admits a Seifert fibration over a hyperbolic 2-orbifold $\Sigma_v$; thus there is a short exact sequence
\[
1 \to \mathbb Z \to \pi_1(M_v) \to \pi_1(\Sigma_v) \to 1
\] where $\mathbb Z$  is the normal cyclic subgroup of $\pi_1(M)$ generated by a fiber. 
If $M_v$ is a hyperbolic piece, then $\pi_1(M_v)$ is hyperbolic relative to  $\{\pi_1(T_1), \ldots, \pi_1(T_\ell) \}$, where $\{T_1, \ldots, T_{\ell} \}$ is the collection of boundary tori of $M_v$. 

Motivated by this structure, Croke and Kleiner introduced the class of \emph{admissible groups} in \cite{CK02}, abstracting the graph of groups structure of graph manifolds. In \cite{MN24}, the authors further introduced the class of \emph{extended admissible groups}, which generalizes fundamental groups of all non-geometric $3$-manifolds as well as Croke--Kleiner admissible groups. In an extended admissible group, vertex groups are allowed to be either central extensions of hyperbolic groups by $\mathbb Z$ or toral relatively hyperbolic groups, yielding a significantly broader and more flexible class of groups.   For the precise definition of extended admissible groups, we refer the reader to Definition~\ref{defn:extended}.

The large-scale geometry of extended admissible groups has recently attracted attention. The main result of \cite{MN24} establishes quasi-isometric rigidity for this class, extending the seminal work of Kapovich--Leeb \cite{KL97} on graph manifold groups. In \cite{Ngu25}, subgroup separability questions for extended admissible groups are studied. Property (QT) has been examined in \cite{NY23}, \cite{HNY23}, while in \cite{HRSS22}, the authors have shown that admissible groups are hierarchically hyperbolic groups. In \cite{ANR24}, the authors demonstrate that admissible groups are $\mathcal{H}$--inaccessible. Additionally, quasi-isometric rigidity is studied in \cite{MN24}, sublinearly Morse boundaries are studied in \cite{NQ24}, and quasi-redirecting boundaries are studied in \cite{NQ25}.

The present paper continues this line of research by examining several coarse-geometric and subgroup-theoretic properties of extended admissible groups, with the goal of extending classical results from non-geometric $3$-manifold groups to this more general setting.

A natural problem accompanying quasi-isometric rigidity is quasi-isometric classification.

\begin{ques}[Quasi-isometric classification] Given a class $\mathcal C$ of finitely generated groups, determine when two elements of $\mathcal C$ are quasi-isometric.
\end{ques}

Each vertex group of an extended admissible group is either a central extension of a hyperbolic group or is relatively
hyperbolic; we call these type $\mathcal S$ and type $\mathcal{H}$ respectively. An extended admissible group $G$ is called an \emph{admissible group} if it has no vertex group of type $\mathcal{H}$. In \cite{MN24}, the authors show that quasi-isometries between extended admissible groups preserve vertex-group types and the quasi-isometry classes of the associated hyperbolic quotients.
This result implies that there are infinitely many quasi-isometry classes of admissible and extended admissible groups. However, it leaves open an important structural question: to what extent does the choice of edge gluing isomorphisms influence the quasi-isometry type? This motivates the following question.

\begin{ques}
    \label{prob:changinggluing}
To what extent do the gluing edge isomorphisms influence the quasi-isometry type
of the resulting admissible groups? 
\end{ques}

Our first result in this paper gives a positive answer to Question~\ref{prob:changinggluing}.

\begin{restatable}{thm}{thm:maintheorem1}
\label{thm:maintheorem1}
  Let $\mathcal{G}$ and $\mathcal{G}'$ be admissible graphs of groups with identical underlying graph, vertex groups, and edge groups, differing only in their edge isomorphisms. Then their fundamental groups are quasi-isometric.
\end{restatable}

This result shows that, within the class of admissible groups, the large-scale geometry is insensitive to the specific gluing data, paralleling classical results for graph manifold groups. 

Motivated by the notion ``flip graph manifolds'' introduced by Kapovich--Leeb \cite{KL98}, we next study admissible groups acting geometrically on Hadamard spaces via so-called CKA actions. In particular, we consider \emph{flip CKA actions}, which arise from specific choices of edge identifications (see Definition~\ref{defn:flipadaction}). These actions play a central role in understanding large-scale nonpositive curvature phenomena.

A key ingredient in this analysis is the notion of omnipotence (see Definition~\ref{defn:omnipotence}), introduced by Wise \cite{Wis00} which has been widely used in subgroup separability. Many familiar classes of groups, including free groups, surface groups, Fuchsian groups, and virtually special hyperbolic groups, are omnipotent. Using omnipotence assumptions on the hyperbolic quotients of vertex groups, we answer a question posed in \cite{NY23}.

\begin{restatable}{thm}{maintheorem}
\label{maintheorem}
       Let $G$ be an admissible group such that each vertex group is a central extension of an omnipotent hyperbolic CAT(0) group by $\mathbb Z$. Then $G$ is quasi-isometric to a group admitting a flip CKA action.
\end{restatable}

As a consequence, we obtain new information about large-scale curvature invariants of extended admissible groups. In particular, we compute their divergence. Divergence is a quasi-isometry invariant introduced by Gersten \cite{Ger94}, which plays a key role in distinguishing geometric behaviors of groups. Gersten showed that non-geometric $3$-manifold groups have either quadratic or exponential divergence, depending on the presence of hyperbolic pieces. We extend this result to extended admissible groups under mild hypotheses on vertex groups.

\begin{restatable}{cor}{cordivergence}
    \label{cordivergence}
     Let $G$ be an extended admissible group such that for each vertex group $G_{v}$ of type $\mathcal{S}$, its non-elementary hyperbolic factor $Q_v$ is omnipotent and is a CAT(0) group. Then the divergence of $G$ is quadratic if and only if $G$ contains no vertex groups of type $\mathcal{H}$ and it is exponential otherwise.

    In particular, assume that $G'$ is another extended admissible group satisfying the same conditions as $G$.  If $G$ contains no vertex groups of type $\mathcal{H}$ and $G'$ contains at least one vertex group of type $\mathcal{H}$, then $G$ and $G'$ are not quasi-isometric.
\end{restatable}

Another application concerns quasi-redirecting boundaries, recently introduced by Qing and Rafi \cite{QR24} as a candidate for a quasi-isometry invariant boundary theory extending the Gromov boundary. While the existence of such boundaries is known for several important classes of groups, it remains open in general. Our main theorem allows the extension of the construction in \cite{NQ25} from  groups admitting  flip CKA actions to a broader class of admissible and extended admissible groups.

\begin{restatable}{cor}{corQRboundary}
 \label{corQRboundary}
      Let $G$ be an extended admissible group such that for each vertex group $G_{v}$ of type $\mathcal{S}$, its non-elementary hyperbolic factor $Q_v$ is omnipotent and is a CAT(0) group. Then $G$ has well-defined quasi-redirecting boundary.
\end{restatable}

In particular, we obtain new examples among free-by-cyclic groups, a class that has been extensively studied but for which the existence of well-defined quasi-redirecting boundaries was previously unknown.

\begin{restatable}{cor}{freebycyclic}
\label{freebycyclic}
Let $\Phi$ be a linearly growing  automorphism of the finite rank free group $F$ and let $G=F\rtimes_\Phi\langle t\rangle$ be its mapping torus. Suppose that $G$ is unbranched in the sense of \cite{BGGH25}. Then $G$ has well-defined quasi-redirecting boundary.
\end{restatable}

Finally, we investigate subgroup structure in extended admissible groups. We study the relationship between strong quasiconvexity (see Definition~\ref{defn:quasiconvex}), finite height , and (virtual) malnormality (see Definition~\ref{defn:malnormalheight}). Building on work of Tran \cite{Tra19} and Hruska--Wise \cite{HW09}, we show that separable, strongly quasiconvex subgroups are virtually almost malnormal in finitely generated groups. As an application, we characterize strongly quasiconvex subgroups of graph manifold groups as precisely those that are virtually malnormal.

\begin{restatable}{thm}{separable}
\label{separable}
Let $H$ be a separable, strongly quasiconvex subgroup of a finitely generated group $G$. Then there is a finite index subgroup $K$ of $G$ containing $H$ such that $H$ is almost malnormal in $K$. Furthermore, suppose that $G$ is virtually torsion-free then $H$ is virtually malnormal. 

Suppose $G$ is the fundamental group of a graph 3-manifold $M$. Then a finitely generated subgroup $H$ of $\pi_1(M)$ is strongly quasiconvex if and only if $H$ is virtually malnormal in $G$.
\end{restatable}

We conclude by establishing several embedding obstructions for extended admissible groups, including consequences of the Rapid Decay property and examples of non-subgroup-separable behavior.

\begin{restatable}{prop}{embeddedsubgroups}
\label{embeddedsubgroups}
Let $G$ be an extended admissible group. Then
\begin{enumerate}
    \item $G$ has Rapid Decay property. As a consequence, amenable groups with exponential growth, Thompson's groups, $SL_{n}(\mathbb Z)$ with $n \ge 3$, intermediate growth groups, Baumslag-Solitar group $BS(p, q)$ (where $p\ne q$) cannot be embedded as subgroups of $G$.

    \item Consider the following \textit{Croke--Kleiner group} $
L = \langle i, j, k, l \mid [i,j], [j,k], [k,l] \rangle$. Suppose that $G$ contains at least one maximal admissible component then there is an embedding $L \to G$ and hence $G$ is not subgroup separable.
\end{enumerate}
\end{restatable}

\subsection{Overview} In Section~\ref{sec:preliminary} we recall the necessary background on trees of spaces and extended admissible groups. Section~\ref{sec:nonpositive} proves the quasi-isometric invariance under changing edge maps. In Section~\ref{sec:CAT(0)} we establish large-scale CAT(0) geometry via flip CKA actions and derive consequences for divergence and boundaries. Section~\ref{sec:subgroups} studies subgroup structure, focusing on strong quasiconvexity and malnormality.

\section{Preliminaries}
\label{sec:preliminary}
In this section, we review some concepts in geometric group theory that will be used throughout the paper.

\subsection{Coarse geometry}
\begin{defn}
    Let $X$ and $Y$ be metric spaces and $f$ be a map from $X$ to $Y$.
    \begin{enumerate}
        \item We say that $f$ is a \emph{$(K,A)$--quasi-isometric embedding} if  for all $x, x'\in X$,
	      \[
		      \frac{1}{K} d(x, x') - A \le d(f(x), f(x')) \le K d(x,x') + A .
	      \]
          \item We say that $f$ is a \emph{$(K,A)$--quasi-isometry} if it is a $(K,A)$--quasi-isometric embedding such that $Y = N_{A}(f(X))$.

          \item We say two quasi-isometries $f, g \colon X \to Y$ are \emph{$A$--close} if
	      \[
		      \sup_{x \in X} d_{Y} (f(x), g(x)) \le A
	      \] and are \emph{close} if they are $A$--close for some $A \ge 0$.
    \end{enumerate}
\end{defn}

\subsection{Tree of spaces of graph of groups} \label{sec:TofS}
We assume familiarity with Bass--Serre theory; see~\cite{SW79} for details. However, to fix notation and terminology, we give some brief definitions.

We first establish some terminology regarding graphs. A \emph{graph} $\Gamma$ consists of a set $V\Gamma$ of vertices, a set $E\Gamma$ of oriented edges, and maps $\iota,\tau:E\Gamma\to V\Gamma$. There is a fixed-point free involution $E\Gamma\to E\Gamma$, taking an edge $e\in E\Gamma$ such that $\iota e=v$ and $\tau e=w$ to an edge $\bar e$ satisfying $\iota \bar e=w$ and $\tau \bar e=v$. We also write $e_+$ and $e_-$ to denote $\tau e$ and $\iota e$ respectively. An \emph{unoriented edge} of $\Gamma$ is the pair $\{e,\bar e\}$. If $v$ is a vertex, we define $\lk(v)=\{e\in E\Gamma\mid e_-=v\}$.

Each connected graph can be identified with a  metric space by equipping its topological realization with the path metric in which each edge has length one. A \emph{combinatorial path} in $X$ is a path $p:[0,n]\to X$ for some $n\in \mathbb N$ such that for every integer $i$, $p(i)$ is a vertex, and $p|_{[i,i+1]}$ is either constant or traverses an edge of $X$ at unit speed. Every geodesic between vertices of $X$ is necessarily a combinatorial path.

\begin{defn}
A \emph{graph of groups} $\mathcal{G} = (\Gamma, \{G_{\hat v}\}, \{G_{\hat e} \}, \{\tau_{\hat e} \})$ consists of the following data:
\begin{enumerate}
    \item a graph $\Gamma$ (called the \emph{underlying graph}),
    \item a group $G_{\hat v}$ for each vertex $\hat{v} \in V (\Gamma)$ (called a \emph{vertex group}),
    \item a subgroup $G_{\hat e} \leq G_{\hat{e}_-}$ for each edge $\hat{e} \in E (\Gamma)$ (called an \emph{edge group}),
    \item an isomorphism $\tau_{\hat{e}} \colon {G}_{\hat{e}} \to {G}_{\bar{\hat{e}}}$  for each $\hat{e} \in E(\Gamma)$ such that $\tau^{-1}_{\hat e}=\tau_{\bar{\hat e}}$ (called an \emph{edge map}).
\end{enumerate}
\end{defn}

The \emph{fundamental group} $G=\pi_{1}(\mathcal G)$ of a graph of groups $\mathcal G$ is as defined in~\cite{SW79}. Via the construction of $G$, we will always view  vertex and edge groups of $\mathcal G$ as subgroups of $G$.

We use the following notation for trees of spaces, similar to~\cite{CM17}.

\begin{defn}
\label{defn:treeofspaces}
\emph{A tree of spaces} $X:=X\left(T,\left\{X_v\right\}_{v \in V (T)},\left\{X_e\right\}_{e \in E (T)},\left\{\alpha_e\right\}_{e \in E (T)}\right)$ consists of:
\begin{enumerate}
   \item a simplicial tree $T$, called the \emph{base tree};
		\item a connected graph $X_v$ for each vertex $v$ of $T$, called a \emph{vertex space};
		\item a connected subgraph $X_e \subseteq X_{e_{-}}$ for each oriented edge $e$ (with the initial vertex denoted by $e_{-}$) of $T$, called an \emph{edge space};
		\item graph isomorphisms  $\alpha_e: X_e \rightarrow X_{\bar{e}}$ for each edge $e \in E T$, such that $\alpha_{\bar{e}}= \alpha_e^{-1}$.
\end{enumerate}
\end{defn}
\textbf{Metrics on $X$:}
 We think of $X$ as a metric space by equipping it with the path metric. Each vertex and edge space $X_x$ of $X$ with $x\in VT\sqcup ET$ is thus endowed with two a priori different metrics: the induced path metric on $X_x$, and the subspace metric when $X_x$ is considered as a subspace of $X$.

\subsection{Tree of spaces from graph of groups} \label{subsec:tree-of-spaces} We now explain how to associate a tree of spaces to a graph of finitely generated groups.

Let $\mathcal{G} = (\Gamma, \{G_{\hat v}\}, \{G_{\hat e} \}, \{\tau_{\hat e} \})$ be a graph of finitely generated groups and let $G$ be the fundamental group of this graph of groups. We recall the associated Bass--Serre tree $T$ is constructed so that vertices (resp.\ edges) of $T$ correspond to left cosets of vertex (resp.\ edge) groups of $\cG$.

We now describe a tree of spaces $X$.  For each $\hat x\in V\Gamma \sqcup E\Gamma$, we fix a finite generating set $J_{\hat x}$ of $G_{\hat x}$, chosen such that $\tau_{\hat e}(J_{\hat e})=J_{\bar {\hat e}}$, and $J_{\hat e}\subseteq J_{\hat v}$ if $\hat e\in ET$ with $\iota {\hat e}= \hat v$.
We now define a graph $W$ with vertex set $V\Gamma\times G$ and edge set \[\{((\hat v, g),(\hat v,gs))\mid g\in G, s\in J_{\hat v}\}.\] The components of $W$ are in bijective correspondence with left cosets of vertex groups of $\cG$, and hence with vertices of $T$. If $v\in VT$ corresponds to $gG_{\hat v}$, we define $X_{v}$  to be the  component of $W$ with vertex set $\{(\hat v,h)\mid h\in gG_{\hat v}\}$. We note that the component of $W$ corresponding to a coset $gG_{\hat v}$ is isometric to the Cayley graph of $G_{\hat v}$ with respect to $J_{\hat v}$.

Suppose $e\in ET$ corresponds to a coset $gG_{\hat e}$. By the definition of $T$, if $\hat v={\hat e}_-$ and $\hat w= {\hat e}_+$, then $v : = e_-$ and $w : = e_+$ correspond to the cosets $gG_{\hat v}$ and $gG_{\hat w}$. We define the edge space $X_{e}$ to  be the graph with vertex set \[\{(\hat v,h)\mid h\in gG_{\hat e}\}\subseteq X_{v}\] and edge set \[\{((\hat v, h),(\hat v,hs))\mid h\in gG_{\hat e}, s\in J_{\hat e}\}.\] Thus $X_e$ is isomorphic to the Cayley graph of $G_{\hat e}$ with respect to $J_{\hat e}$.  The attaching map $\alpha_{e}:X_{e}\to X_{w}$ is defined by  $\alpha_{e}:(v,h)\mapsto (w,g\tau_{\hat e}({g^{-1}h}))$ on vertices, and similarly on edges, where $\tau_{\hat e}:G_{\hat e}\to G_{\bar{\hat e}}\leq G_{\hat w}$ is the edge map of $\cG$. 
\begin{defn}
    \label{defn:treeofspacesassociated}
	Given a graph of finitely generated groups $\cG$, the tree of spaces $X$ constructed above  is \emph{the tree of spaces associated with the graph of groups $\cG$}.
\end{defn}

The tree of spaces $X$ is a proper geodesic metric space (see Lemma 2.13 of~\cite{CM17}). The natural action of $G$ on $W$ (fixing the $V\Gamma$ factor) induces an action of $G$ on $X$. Applying the  Milnor--Schwarz lemma we deduce:

\begin{prop} 
[Section 2.5 of~\cite{CM17}]\label{prop:treeofspacesBST}
	Suppose $G$, $T$ and $X$ are as above. Then there exists a quasi-isometry $f: G \rightarrow X$ and $A\geq 0$ such that
	$d_{\mathrm{Haus }}\left(f\left(gG_{\hat x}\right), X_{x}\right) \leq A$  for all $x \in V T\sqcup  E T$, where $x$ corresponds to the coset $gG_{\hat x}$.
\end{prop}

The following theorem explains how to build a quasi-isometry between trees of spaces
by patching together quasi-isometries of vertex spaces. This can be done if quasi-isometries
on adjacent vertex spaces agree up to a uniformly bounded error on their common edge
space.

\begin{thm}\cite[Corollary~2.16]{CM17}
\label{thm:CM17}
Let $K \geq 1$ and $A \geq 0$. Suppose that $X:= X \bigl (T,\left\{X_v\right\},\left\{X_e\right\},\left\{\alpha_e\right\} \bigr )$ and $X' :=X' \bigl (T',\left\{X_v^{\prime}\right\},\left\{X_e^{\prime}\right\},\left\{\alpha_e^{\prime}\right\} \bigr )$ are trees of spaces, and that there is a tree isomorphism $\xi \colon  T \rightarrow T^{\prime}$. Suppose for every $v \in V (T)$ and $e \in E (T)$ there is a $(K, A)$--quasi-isometry $\phi_v \colon X_v \to X_{\xi(v)}^{\prime}$ and $\phi_e \colon X_e \to X_{\xi(e)}^{\prime}$. Suppose also that for every $e \in E (T)$, the following diagrams commute up to uniformly bounded error $A$. 
\[\begin{tikzcd}
		{X_e} & {X'_{\xi(e)}} \\
		{X_{e_{-}}} & {X'_{\xi(e_{-})}} 
		\arrow["\phi_e", from=1-1, to=1-2]
		\arrow[, from=1-1, to=2-1]
		\arrow["\phi_{e_{-}}"', from=2-1, to=2-2]
		\arrow[, from=1-2, to=2-2]
	\end{tikzcd}\qquad 
  \begin{tikzcd}
		{X_e} & {X'_{\xi(e)}} \\
		{X_{e_{+}}} & {X'_{\xi(e_{+})}} 
		\arrow["\phi_e", from=1-1, to=1-2]
		\arrow["\alpha_{e}", from=1-1, to=2-1]
		\arrow["\phi_{e_{+}}"', from=2-1, to=2-2]
		\arrow["\alpha'_{\xi(e)}", from=1-2, to=2-2]
	\end{tikzcd}\]
Then there is a quasi-isometry $\phi: X \rightarrow X^{\prime}$ such that $\phi |_{X_v}=\phi_v$ for every $v \in V (T)$.
\end{thm}

\subsection{Extended admissible groups}
We now define the class of extended admissible groups.

\begin{defn}
\label{defn:extended}
   A group $G$ is an  \emph{extended admissible group} if it is the fundamental group of a graph of groups $\cG$ such that:
	\begin{enumerate}
		\item The underlying graph $\Gamma$ of $\cG$ is a connected finite graph with at least one edge, and every edge group is $\Z^2$.
		\item Each vertex group $G_v$ is one of the following two types:
		      \begin{enumerate}
			      \item  Type $\mathcal{S}$: $G_v$ has center $Z_v : = Z(G_v) \cong \mathbb Z$  such that the quotient $Q_v := G_v / Z_v$ is a non-elementary hyperbolic group. We call $Z_v$ and $Q_v$  the \emph{kernel} and \emph{hyperbolic quotient} of $G_v$ respectively.
			      \item Type $\mathcal{H}$: $G_v$ is hyperbolic relative to a collection $\PP_v$ of virtually $\Z^2$-subgroups, where all edge groups incident to $G_v$ are contained in $\PP_v$, and  $G_v$ doesn't split relative to $\PP_v$  over a subgroup of an element of $\PP_v$.
		      \end{enumerate}
		\item  For each vertex group $G_v$, if $e,e'\in \lk(v)$ and $g\in G_v$, then  $gG_eg^{-1}$ is commensurable to $G_{e'}$ if and only if both $e=e'$ and $g\in G_e$.

		\item For every edge group ${ G}_e$ such that $G_{e_{-}}$ and $G_{e_{+}}$ are vertex groups of type $\mathcal{S}$, the subgroup generated by $\tau_{\bar e}(Z_{e_+}\cap {G}_{\bar e})$  and $Z_{e_-}\cap G_e$ has finite index in ${ G}_e$.
	\end{enumerate}
 \end{defn}

\begin{defn}
    \label{defn:admissible}
	An extended admissible group $G$ is called an \emph{admissible group} if it has no vertex group of type $\mathcal{H}$.
\end{defn}
 
\textbf{Convention:}
	For the rest of this paper, if $G$ is an extended admissible group, we will assume that all the data $\mathcal{G}$, $G_{\hat v}$, $Z_{\hat v}$, $Q_{\hat v}$, etc. in Definition \ref{defn:extended} are fixed, and will make use of this notation without explanation.  If $G'$ is another admissible group, we use the notation  $\mathcal{G}'$, $G'_{\hat v}$, $Z'_{\hat v}$, $Q'_{\hat v}$ etc.

Below are some examples of extended admissible groups. 

\begin{exmp}
\label{exmpleextended}
    \begin{enumerate}
        \item (3-manifold groups) The fundamental group of a compact, orientable, non-geometric, irreducible 3-manifold $M$ with
		      empty or toroidal boundary is an extended admissible group. Seifert fibered and hyperbolic pieces correspond to type $\cS$ and $\cH$ vertex respectively. Fundamental groups of graph manifolds are admissible groups.
\item (Torus complexes) 
Let $n \ge 3$ be an integer. Let $T_1, T_2, \ldots, T_n$ be a family of flat two-dimensional tori. For each $i$, we choose a pair of simple closed geodesics $a_i$ and $b_i$ such that $a_i \cap b_i \neq \emptyset$ and $\operatorname{length} (b_i) = \operatorname{length} (a_{i+1})$, identifying $b_i$ and $a_{i+1}$ and denote the resulting space by $X$. For each $i \in \{1, \ldots, n-1\}$, we denote $V_{i} : = T_{i} \cup T_{i+1} / \{ b_i =  a_{i+1} \}$. Let $S^{1}_{i} \subset V_i$ be the subspace of $V_i$ obtained by gluing $b_i$ to $a_{i+1}$. The space $X$ is obtained by gluing each $V_i$ to $V_{i+1}$ via the gluing map $$\tau_{i} \colon b_{i+1} \times S^{1}_i \subset V_i \to a_{i+1} \times S^{1}_{i+1} \subset V_{i+1}$$ by sending $b_{i+1} \to S^{1}_{i+1}$ and $S^{1}_{i}\to a_{i+1}$ accordingly. Such a gluing map is called a ``flip'' map in the literature.

Note that $V_i$ is homotopic equivalent to the product of $S^{1}_i$ with the wedge of two circles $a_i$ and $b_{i+1}$. 
The fundamental group $G = \pi_{1}(X)$ has a graph of groups structure where each vertex group $\pi_1(V_i) = (\langle a_i \rangle * \langle b_{i+1} \rangle ) \times \Z = F_{2} \times \Z$, edge groups are $\Z^2$ and edge maps are induced by the gluing maps $\tau_i$.
It is clear that with this graph of groups structure, $\pi_1(X)$ is an admissible group.

Note that our space $X$ is a local CAT(0) space, and hence the universal cover $\widetilde{X}$ is CAT(0) by the Cartan-Hadamard theorem. This space is studied in \cite{CK00}.
    \end{enumerate}
\end{exmp}

\begin{lem} \cite[Lemma~4.2]{HRSS22}
    \label{lem:HRSS22}
Let $\mathcal{G} = (\Gamma, \{G_{\hat{v}}\}, \{G_{\hat{e}} \}, \{\tau_{\hat{e}} \})$ be an admissible group. Each vertex group $G_{\hat v}$  has an infinite generating set $S_{\hat v}$ so that
the following holds.
\begin{enumerate}
    \item The Cayley graph $\mathrm{Cay}(G_{\hat{v}}, S_{\hat{v}})$ is quasi-isometric to a line.

    \item The inclusion map $Z_{\hat{v}} \to \mathrm{Cay}(G_{\hat{v}}, S_{\hat{v}})$ is a $Z_{\hat{v}}$--equivariant quasi-isometry.
\end{enumerate}
\end{lem}

\begin{rem}
     Without loss of generality, we can assume that the finite generating set $J_{\hat{v}}$ of $G_{\hat{v}}$ is contained in $S_{\hat{v}}$.
\end{rem}

Recall from the construction in Section \ref{sec:TofS} that each vertex space $X_v$ of $X$ is identified with the Cayley graph of a vertex group $G_{\hat v}$ of $\mathcal{G}$ with respect to some generating set $J_{\hat v}$.

\begin{defn} (Subspace $L_v$ and $\mathcal{H}_v$)
\label{defn:subspaceLv}
 Suppose that $v \in T$ corresponds to a coset $gG_{\hat{v}}$. Let
$L_v \subset X_v$
be the graph with vertex set $gG_{\hat{v}}$ and with an edge connecting $x, y\in gG_{\hat{v}}$ if $x^{-1}y\in S_{\hat{v}}$.  In particular, $L_v$ is isometric to $\mathrm{Cay}(G_{\hat{v}}, S_{\hat{v}})$, which is a quasi-line by Lemma~\ref{lem:HRSS22}. %The corresponding subspace $L'_v$ is defined  analogously in $X'_v$.    

Let $\mathcal{H}_v$ be the graph with vertex set $g G_{\hat v}$ and an edge connecting  $x, y \in g G_{\hat v}$ if $x^{-1} y \in J_{\hat v} \cup Z_{\hat v}$. It is isometric to $\mathrm{Cay}(G_{\hat v}, J_{\hat v} \cup Z_{\hat v})$. We call $\mathcal{H}_v$ is the \textit{quotient space} of $X_v$.
\end{defn}

\begin{rem}
    \label{rem:equivariant}
    For any $g \in G$ and for each vertex $v \in V(T)$ we have $g L_v = L_{g v}$.
\end{rem}

\begin{defn} [Quotient maps, boundary lines]
\label{defn:somemaps}
Suppose that $v \in T$ corresponds to a coset $gG_{\hat{v}}$. Since $L_v$ and $\mathcal H_v$ are each obtained from $X_v$ by adding extra edges, there are distance non-increasing maps $p_v \colon X_v \to L_v$ and $\pi_v\colon X_v\to \mathcal H_v$ that are the identity on vertices.
 We call such $\pi_v:X_v\to \mathcal{H}_v$ is a \textit{quotient map}. For each $e\in E(T)$ with $v=e_-$, we define the  \textit{boundary line} $\ell_e$ of $\mathcal{H}_v$ associated to $e$ is $$\ell_e:= \pi_v(X_e)\subseteq \mathcal{H}_v.$$
 Let $w$ be an adjacent vertex of $v$ and denote the oriented edge $[v,w]$ by $e$. Let $$\psi_e \colon \ell_{\bar{e}} \to L_v$$ be the restriction to the boundary line $\ell_{\bar{e}}$  of the composition
$
p_{v} \circ  \alpha_{\bar{e}} \circ \pi_{w}^{-1}
$. 
\end{defn}

\begin{rem}
\begin{enumerate}
    \item   The space $\mathcal H_v$ is constructed to represent the geometry of $Q_{\hat{v}} = G_{\hat{v}}/Z_{\hat{v}}$ and is relatively hyperbolic  to the collection $\{ \ell_e \}_{e_{-} = v}$ (see \cite[Lemma~2.15]{HRSS22}).

    \item It is proved in \cite[Lemma~2.18]{ANR24} that $\psi_e$ is a uniform quasi-isometry. Namely, there exists constants $\lambda \ge 1, c \ge 0$ such that for each oriented edge $e$ in $T$ then $\psi_e \colon \ell_{\bar{e}} \to L_v$ is a $(\lambda, c)$--quasi-isometry.
\end{enumerate}
\end{rem}

\begin{lem}
\label{lem:uniformQI}
    There exist constants $\lambda \ge 1, c \ge 0$ such that the following holds. Suppose that $v \in T$ corresponds to a coset $gG_{\hat{v}}$ where $\hat v$ is a vertex in the underlying graph $\Gamma$. Consider the map $$f_{v} \colon X_v \to \mathcal{H}_v \times L_v$$ defined by $x \mapsto (\pi_{v}(x), p_{v} (x))$ where $\pi_v$ and $p_v$ are maps given by Definition~\ref{defn:somemaps}. Then $f_v$ is a $(\lambda, c)$--quasi-isometry.
\end{lem}

\begin{proof}
 We consider  two natural actions $G_{\hat{v}} \curvearrowright Q_{\hat{v}}$ and $G_{\hat{v}} \curvearrowright L_{\hat{v}} := \mathrm{Cay} (G_{\hat{v}}, S_{\hat{v}})$ of $G_{\hat{v}}$ on quotients $Q_{\hat{v}}$ and the quasi-line $\mathrm{Cay} (G_{\hat{v}}, S_{\hat{v}})$ respectively.
It is shown in \cite[Corollary~4.3]{HRSS22} that the diagonal action $G_{\hat{v}} \curvearrowright Q_{\hat{v}} \times L_{\hat{v}}$  is metrically proper and co-bounded, and hence  the orbit map (with respect to a fixed basepoint) denoted by 
$
f_{\hat{v}} \colon G_{\hat{v}} \to Q_{\hat{v}} \times L_{\hat{v}}
$
is a quasi-isometry such that the composition of $f_{\hat{v}}$ with the projection $Q_{\hat{v}} \times L_{\hat{v}} \to Q_{\hat{v}}$ is the quotient map $q_{\hat{v}} \colon G_{\hat{v}} \to Q_{\hat{v}} = G_{\hat{v}} / Z_{\hat v}$.
It implies that  $f_v$ is a quasi-isometry.
Since there are finitely many vertices in the underlying graph $\Gamma$, we conclude that $f_v$ is a quasi-isometry with uniform quasi-isometric constants $\lambda, c$.
\end{proof}

\section{Changing edge maps does not change quasi-isometric type}
\label{sec:nonpositive}

In this section, we are going to prove Theorem~\ref{thm:maintheorem1} by showing that if two admissible groups $\mathcal{G} = (\Gamma, \{G_{\hat{v}}\}, \{G_{\hat{e}} \}, \{\tau_{\hat{e}} \})$ and $\mathcal{G}' = (\Gamma, \{G_{\hat{v}}\}, \{G_{\hat{e}} \}, \{\sigma_{\hat{e}} \})$ differ only in their edge isomorphisms then $G = \pi_1(\mathcal{G})$ and $G' = \pi_1(\mathcal{G}')$ are quasi-isometric. 
Fix trees of spaces $(X,T)$, $(X',T)$ associated with admissible groups $G$ and $G'$ respectively, with the same associated Bass--Serre tree $T$. By Proposition~\ref{prop:treeofspacesBST}, $G$ and $G'$ are quasi-isometric to $X$ and $X'$ respectively. Hence it suffices to show that $X$ and $X'$ are quasi-isometric. To do so, we are going to construct collections of quasi-isometries 
\[
\{\phi_v : X_v \to X'_v\}_{v\in V(T)} \quad \text{and} \quad \{\phi_e : X_e \to X'_e\}_{e\in E(T)}
\]
between the vertex and edge spaces in tree of spaces so that these collections of maps satisfy conditions Theorem~\ref{thm:CM17}.

%The overall strategy is as follows:
%\begin{enumerate}
 %   \item \textbf{Useful Maps:} We first introduce canonical maps on certain vertex spaces and their ``boundary lines''. These maps are quasi-isometries with uniform quasi-isometric constants.
    
 %   \item \textbf{Construction of Vertex/Edge Quasi-Isometries:} We then use these maps (and a fixed choice of quasi-isometries on representative vertex orbits) to extend to quasi-isometries on all vertex spaces and, by restriction, on the edge spaces.
 %   \item \textbf{Verification of Commutativity:} Finally, we check that the maps satisfy the compatibility conditions in Theorem~\ref{thm:CM17}. This guarantees the existence of a global quasi-isometry \(\varphi: X \to X'\) and, hence, the quasi-isometry between \(G\) and \(G'\).
%\end{enumerate}

\textbf{Outline of the proof:} In Section~\ref{subsec:construction}, using a fixed choice of quasi-isometries on representatives of vertex orbits, we extend these maps equivariantly to all vertex spaces and, by restriction, obtain induced quasi-isometries on edge spaces. The main technical step is to verify that the resulting vertex and edge quasi-isometries satisfy the compatibility conditions of Theorem~\ref{thm:CM17}, despite the change in edge isomorphisms. This is done in Section~\ref{sec:compatible}.

\subsection{Construction of vertex/edge maps}
\label{subsec:construction}

Since \(G\) and \(G'\) share the same Bass–Serre tree \(T\), their vertex spaces \(X_v\) and \(X'_v\) (and similarly, the edge spaces \(X_e\) and \(X'_e\)) are naturally identified. The only difference is in the gluing isomorphisms (from \(\tau_e\) to \(\sigma_e\)).
 Let $G_{\hat{v}_0}, \ldots, G_{\hat{v}_m}$ be the vertex subgroups of $G$. We proceed as follows:

     \underline{Choice of Transversals:} For each $i$,  we fix $\mathfrak S_{\hat{v}_i}$ a set of transervals for left cosets of $G_{\hat{v}_i}$ in $G$ such that $1 \in \mathfrak S_{\hat{v}_i}$.

    \underline{Base quasi-isometries:}  For each vertex $v$ in $V(T)$, let $L_v$ and $L'_v$ be the spaces defined in Definition~\ref{defn:subspaceLv} with respect to tree of spaces $X$ and $X'$. Given a vertex $\hat{v}_i \in \{\hat{v}_0, \hat{v}_1, \ldots, \hat{v}_m \}$ in the underlying graph $\Gamma$, we fix a vertex $v_i \in V(T)$ corresponding to $G_{\hat{v}_i} = 1 \cdot G_{\hat{v}_i}$. Recall from Definition~\ref{defn:subspaceLv} that $L_{v_i}$ (resp $L'_{v_i}$) is the graph with vertex set $1 \cdot G_{\hat{v}_i}$ with an edge connecting $x, y \in 1\cdot G_{\hat{v}_i}$ if $x^{-1} y \in S_{\hat{v}_i}$. In particular, $L_{v_i} = L'_{v_i}$. We also recall two quotient spaces $\mathcal{H}_{v_i}$ and $\mathcal{H}'_{v_i}$ from Definition~\ref{defn:subspaceLv} as well and remark that $\mathcal{H}_{v_i} = \mathcal{H}'_{v_i}$.
    We fix a quasi-isometry $\zeta_{v_i} \colon L_{v_i} \to L'_{v_i}$ which is the composition: 
$$
 L_{v_i} \to X_{v_i}  \to \mathcal{H}_{v_i} \times L_{v_i} = \mathcal{H}'_{v_i} \times L'_{v_i} \to L'_{v_i}
$$ where the first map is the inclusion of $L_{v_i}$ to $X_{v_i}$ (as $L_v \subset X_v$), the second map is $f_{v_i}$ given by Lemma~\ref{lem:uniformQI} and the third map is the natural projection of $\mathcal{H}'_{v_i} \times L'_{v_i}$ into its second factor.

 \underline{Extending to all vertices:} 
At the moment, we have defined maps $\zeta_{v_0}, \zeta_{v_1}, \ldots, \zeta_{v_m}$. We need to define $\zeta_{v}$ for an arbitrary vertex $v$ in $V(T)$. For each vertex $v \in V(T)$, there exists a vertex $\hat{v}_i \in \{ \hat{v}_0, \hat{v}_1, \ldots, \hat{v}_m \}$ and a group element $t \in \mathfrak S_{\hat{v}_i}$ such that $v$ corresponds to the coset $t G_{\hat{v}_i}$.

We write $t$ in reduced form relative to a fixed maximal tree $\Lambda \subset \Gamma$. Namely, fix a maximal tree $\Lambda \subset \Gamma$.
$G$ has a finite generating set of the form $\mathcal{S} = \cup_{i=1}^{m} J_{\hat{v}_i} \cup J_0$ where $J_0$ consists of stable letters $t_e$ corresponding to egdes outside the maximal tree $\Lambda$ (and $t_e = 1$ when $e \in E(\Lambda)$). Similarly for $G'$ with the same maximal tree $\Lambda$.
  We first write the group element $t \in \mathfrak S_{\hat{v}_i}$ in reduced form $$t = g_{0} t_{\alpha_1} g_{1} t_{\alpha_2} \ldots t_{\alpha_k} g_k$$ where each $g_i$ is a group element in a vertex group of $G$ and $\alpha_1 \cdots \alpha_k$ is a loop in $\Gamma$.  We then define $$t' : =g_{0} t'_{\alpha_1} g_{1} t'_{\alpha_2} \ldots t'_{\alpha_k} g_k$$ which is a group element in $G'$. Note that $t v_i = t' v_i = v$ in the Bass-Serre tree $T$ and hence $L'_v = L'_{t'v_i} = t' L'_{v_i}$ and $L_v = L_{t v_i} = t L_{v_i}$ by Remark~\ref{rem:equivariant}. Since $L_v = t L_{v_i}$, it follows that each element in $L_v$ can be written as $tx$ for some $x \in L_{v_i}$. This yields a well-defined quasi-isometry 
  \[
    \zeta_v: L_v \to L'_v,\quad t x \mapsto t'\,\zeta_{v_i}(x).
    \]
    By finiteness of the $G$--orbits of vertices, the map $\zeta_v$ can be chosen uniformly quasi-isometric.

\begin{defn} 
\label{defn:vertex-edge maps}For each vertex $v$ in $V(T)$, let $f_v \colon X_v \to \mathcal{H}_v \times L_v$ and $f'_v \colon X'_v \to \mathcal{H}'_v \times L'_{v}$ be the maps given by Lemma~\ref{lem:uniformQI}. Let $e$ be an edge in $E(T)$ with $e_{-} = v$. We define the \textit{vertex map} \[
 \phi_v: X_v \to X'_v \quad \text{by} \quad \phi_v = (f'_v)^{-1} \circ (\mathrm{id} \times \zeta_v) \circ f_v.
\]
 and  the \textit{ edge map}
\[
\phi_e: X_e \to X'_e \quad \text{by} \quad \phi_e = \Bigl((f'_v)^{-1}\Bigr)\Big|_{X'_e} \circ (\mathrm{id} \times \zeta_v) \circ \Bigl(f_v|_{X_e}\Bigr)
\]
\end{defn}

\subsection{Proof of Theorem~\ref{thm:maintheorem1}}
\label{sec:compatible}
In this section, we are going to prove Theorem~\ref{thm:maintheorem1}. Given two admissible groups $\mathcal{G} = (\Gamma, \{G_{\hat{v}}\}, \{G_{\hat{e}} \}, \{\tau_{\hat{e}} \})$ and $\mathcal{G}' = (\Gamma, \{G_{\hat{v}}\}, \{G_{\hat{e}} \}, \{\sigma_{\hat{e}} \})$. Let $X:= X \bigl (T,\left\{X_v\right\},\left\{X_e\right\},\left\{\alpha_e\right\} \bigr )$ and $X' :=X' \bigl (T',\left\{X_v^{\prime}\right\},\left\{X_e^{\prime}\right\},\left\{\alpha_e^{\prime}\right\} \bigr )$ be the tree of spaces associated to $\mathcal{G}$, $\mathcal{G}'$ given by Section~\ref{subsec:tree-of-spaces}. Here $\alpha_e \colon X_e \to X_{\bar e}$ and $\alpha'_e \colon X'_e \to X'_{\bar e}$.

Let $\{\phi_v\}$ and $\{\phi_e\}$ be the collection of vertex maps and edge maps given by Definition~\ref{defn:vertex-edge maps}.
We fix uniform constants $K\ge 1, A \ge 0$ so that each $\phi_v$ and $\phi_e$ is a $(K,A)$--quasi-isometry. The isomorphism $\xi \colon T \to T$ here we are using is the identity $T \to T$.

For each vertex $v$ in $V(T)$, let $e$ be an edge such that $e_{-} = v$.
By the construction of $\phi_v$ and $\phi_e$ in Section~\ref{subsec:construction}, the following diagram is commuted up to uniformly bounded error
\[\begin{tikzcd}
		{X_e} & {X'_{e}} \\
		{X_{v}} & {X'_{v}} 
		\arrow["\phi_e", from=1-1, to=1-2]
		\arrow["i_e", from=1-1, to=2-1]
		\arrow["\phi_{v}"', from=2-1, to=2-2]
		\arrow["i'_e", from=1-2, to=2-2]
	\end{tikzcd}\]
Here $i_e, i'_e$ are inclusion maps from edge spaces to vertex spaces. Hence we establish commutativity (up to uniformly bounded error) of the first sub-diagram in Theorem~\ref{thm:CM17}.

For the rest of the proof, we are going to verify that our maps satisfy the commutativity (up to uniformly bounded error) of the second sub-diagram in Theorem~\ref{thm:CM17}. In other words, if $w : = e_{+}$ then we verify that the following diagram is commuted up to uniformly bounded error. 
\[\begin{tikzcd}
		{X_e} & {X'_{e}} \\
		{X_{v}} & {X'_{w}} 
		\arrow["\phi_e", from=1-1, to=1-2]
		\arrow["\alpha_{e}", from=1-1, to=2-1]
		\arrow["\phi_{w}"', from=2-1, to=2-2]
		\arrow["\alpha'_{e}", from=1-2, to=2-2]
	\end{tikzcd}\]
\textbf{Claim~1:} There exists a uniform constant  $C > 0$ such that
\begin{enumerate}
    \item For each oriented edge $e = [v,w]$, $X_e$ is quasi-isometric to $L_w \times L_v$ via the following $(C, C)$--quasi-isometry $\rho_e \colon X_e \to L_w \times L_v$ defined by
\[
    x \mapsto ( \psi_{\bar{e}} \circ \pi_{v} (x), \psi_{e} \circ \pi_{w} \circ \alpha_{e} (x)). 
\]
    \item The following diagram is commuted up to an error $C$. \[
    \begin{tikzcd}
		{X_e} & {X_{\bar e}} \\
		{L_{w} \times L_v} & {L_{v} \times L_w} 
		\arrow["\alpha_e", from=1-1, to=1-2]
		\arrow["\rho_e", from=1-1, to=2-1]
		\arrow["\mathfrak{flip}_{(v,w)}", from=2-1, to=2-2]
		\arrow["\rho_{\bar{e}} ", from=1-2, to=2-2]
	\end{tikzcd}
\]  
\end{enumerate}
It is clear that $(2)$ follows from $(1)$. For $(1)$, condition $(4)$ of Definition~\ref{defn:extended} gives us $$X_e\simeq_{q.i} G_e\simeq_{q.i} \langle Z_v,\alpha_{\overline{e}}(Z_w)\rangle = \alpha_{\overline{e}}(Z_w)\times Z_v.$$ Also, we can rewrite $\rho_e$ as $$x\mapsto (\psi_{\overline{e}}\circ \pi_v(x),\psi_e\circ \pi_w\circ \alpha_e(x))=(p_w\circ \alpha_{\overline{e}}(x),p_v(x)),$$ where $p_v:X_v\rightarrow L_v,p_w:X_w\rightarrow L_w.$ Since $Z_v\rightarrow L_v$ is a $Z_v-$equivariant quasi-isometry, the following map is a quasi-isometry 
$$X_e\simeq_{q.i.} \alpha_{\overline{e}}(Z_w) \times Z_v \xrightarrow{p_w\times p_v}L_w\times L_v.$$ 

\textbf{Claim~2:} The following diagram commutes up to a uniform error.
\[\begin{tikzcd}
		{X_e} & {X'_{e}} \\
		{X_{w}} & {X'_{w}} 
		\arrow["\phi_e", from=1-1, to=1-2]
		\arrow["i_{\bar e} \circ \alpha_{e}", from=1-1, to=2-1]
		\arrow["\phi_{w}"', from=2-1, to=2-2]
		\arrow["i'_{\bar e} \circ \alpha'_{e}", from=1-2, to=2-2]
	\end{tikzcd}\]
 For our notations purpose, we write $f  \approx g$ to mean two maps $f$ and $g$ are uniform close. 

 According to the diagram above Claim~1, we have
 \[
 i'_{\bar e} \circ \phi_{\bar{e}} \circ \alpha_{e} \approx \phi_{w} \circ i_{\bar{e}} \circ \alpha_e
 \]
and hence to show that the above diagram is commuted up to a uniform error, it suffices to verify that 
\begin{equation}
\label{equ:chasing}
    \phi_{\bar e} \circ \alpha_e  \approx \alpha'_{e} \circ \phi_e
\end{equation}

To see this, we consider the following diagram: 

\[\begin{tikzcd}
		{X_v} & {X_e} & {X_{\bar{e}}} &{X_w} \\
		{Y_v \times L_v} & {L_w \times L_v} & {L_v \times L_w} & {Y_w \times L_w} \\
        {Y_v \times L'_v} & {L'_w \times L'_v} & {L'_v \times L'_w} & {Y_w \times L'_w} \\
        {X'_v} & {X'_e} & {X'_{\bar{e}}} & {X'_w} \\
		\arrow["i_e", from=1-2, to=1-1]
  \arrow["\alpha_e", from=1-2, to=1-3]
  \arrow["i_{\bar{e}}", from=1-3, to=1-4]
  \arrow["", from=2-2, to=2-3]
  \arrow["", from=2-3, to=2-4]
  \arrow["", from=3-2, to=3-1]
  \arrow["", from=3-2, to=3-3]
  \arrow["", from=3-3, to=3-4]
  \arrow["", from=4-2, to=4-1]
  \arrow["", from=4-2, to=4-3]
  \arrow["", from=4-3, to=4-4]
  \arrow["", from=2-1, to=3-1]
  \arrow["", from=3-1, to=4-1]
  \arrow["", from=2-2, to=3-2]
  \arrow["", from=1-3, to=2-3]
  \arrow["", from=2-3, to=3-3]
  \arrow["", from=1-4, to=2-4]
  \arrow["", from=2-4, to=3-4]
   \arrow["", from=2-2, to=2-1]
  \arrow["", from=3-4, to=4-4]
  \arrow["f_w", from=1-4, to=2-4]
  \arrow["(f'_w)^{-1}", from=3-4, to=4-4]
  \arrow["\rho_{\bar{e}}", from=1-3, to=2-3]
    \arrow["\mathfrak{flip}_{(v,w)}", from=2-2, to=2-3]
    \arrow["\mathfrak{flip'}_{(v,w)}", from=3-2, to=3-3]
    \arrow["f_v", from=1-1, to=2-1]
    \arrow["(f'_v)^{-1}", from=3-1, to=4-1]
\arrow["\rho_e", from=1-2, to=2-2]
\arrow["(\rho'_e)^{-1}", from=3-2, to=4-2]
\arrow["(\rho'_{\bar{e}})^{-1}", from=3-3, to=4-3]
\arrow["i'_{\bar{e}}", from=4-3, to=4-4]
\arrow["i'_e", from=4-2, to=4-1]
\arrow["\zeta_{v} \times \zeta_{w}", from=2-3, to=3-3]
\arrow["\zeta_{w} \times \zeta_{v}", from=2-2, to=3-2]
\arrow["\alpha'_e", from=4-2, to=4-3]
\arrow["\psi_{\bar{e}}^{-1} \times id", from=2-2, to=2-1]
\arrow["\psi_{e}^{-1} \times id", from=2-3, to=2-4]
\arrow["id \times \zeta_v", from=2-1, to=3-1]
\arrow["id \times \zeta_w", from=2-4, to=3-4]
	\end{tikzcd}\]
We note that:
\begin{enumerate}
    \item the compositions of maps in the first and the fourth columns are $\phi_v$ and $\phi_w$ respectively;
    \item the compositions of maps in the second and the third columns are uniformly close to $\phi_e$ and $\phi_{\bar{e}}$ respectively;
    \item  construction of maps $\phi_e$ and $\phi_v$, together with Claim 1, shows that the sub-diagram in the diagram above either commutes or commutes up to a uniform error.
\end{enumerate}
Therefore it is routine to chase around the above diagram to check that $\phi_{\bar e} \circ \alpha_e  \approx \alpha'_{e} \circ \phi_e$, establishing (\ref{equ:chasing}). 
Claim~2 is confirmed.

In conclusion, the collections $\{\phi_v\}_{v \in V(T)}$ and $\{\phi_e\}_{e \in E(T)}$ satisfy the hypotheses of Theorem~\ref{thm:CM17}. Therefore  there is a quasi-isometry $$\phi: X \rightarrow X'$$ such that $\phi |_{X_v}=\phi_v$ for every $v \in V (T)$. %Applying Proposition~\ref{prop:treeofspacesBST}, we  obtain the following result.

\section{Admissible groups are CAT(0) on the large scale}
\label{sec:CAT(0)}
In this section, we use Theorem~\ref{thm:maintheorem1} to  prove Theorem~\ref{maintheorem}.

\subsection{Flip CKA action}
We refer the reader to \cite{CK02}, \cite{NY23} for the material recalled here.
\begin{defn}
    \label{defn:admissibleaction}
We say that the action $G \curvearrowright X$ is  {\it Croke-Kleiner admissible} (CKA) if $G$ is an admissible group, and $X$ is a Hadamard space ({{i.e, a complete proper CAT(0) space}}), and the action is   geometric (i.e., properly and cocompactly by isometries). The space $X$ is called the \emph{admissible space} for the CKA action $G \curvearrowright X$.  
\end{defn}

Let $G \curvearrowright X$ be a {Croke-Kleiner} admissible action, {{where $G$ is the fundamental group of an admissible graph of groups $\mathcal{G}$ and let $G \curvearrowright T$ be the action of $G$ on the associated Bass-Serre tree {{of $\mathcal{G}$}} ({{we refer the reader to Section~2.5 in \cite{CK02} for a brief discussion}}). Let $T^0 = \mathrm{Vertex}(T)$ and $T^1 = \mathrm{Edge}(T)$ be the vertex and edge sets of $T$. For each $\sigma \in T^0 \cup T^1$, we let $G_{\sigma} \le G$ be the stabilizer of $\sigma$. For each vertex $v \in T^0$, let $Y_{v} := \mathrm{Minset}(Z(G_v)) := \cap_{g \in Z(G_v)} \mathrm{Minset}(g)$ and for every edge $e \in E$ we let $Y_{e} := \mathrm{Minset}(Z(G_e)) := \cap_{g \in Z(G_e)} \mathrm{Minset}(g)$.
We note that the assignments $v \to Y_v$ and $e \to Y_e$ are $G$--equivariant with respect to the natural $G$ actions.

%{{By Flat Torus Theorem (for example, see Theorem~II.7.1 in \cite{BH99}), we have:}}
% \begin{lem}
%If $H=\mathbb Z^k$ for some $k\ge 1$ then $\mathrm{Minset}(H)=\cap_{h\in H} \mathrm{Minset}(h)$ splits isometrically as a metric product $C\times \mathbb E^k$ so that $H$ acts trivially on $C$ and as a translation lattice on $\mathbb E^k$. Moreover, $Z(H,G)$ acts cocompactly on $C\times \mathbb E^k$.
%\end{lem}
%As a corollary, we have  

We recall some facts from \cite[Section~3.2]{CK02} and \cite[Section 2]{NY23}.
\begin{enumerate}
    \item $G_v$ acts co-compactly on $Y_v=\bar Y_v\times \mathbb R$ and $Z(G_v)$ acts by translation on the $\mathbb R$--factor and trivially on $\bar Y_v$ where $\bar Y_v$ is a Hadamard space.
    \item $G_e=\mathbb Z^2$ acts co-compactly on $Y_e=\bar Y_e \times \mathbb R^2\subset Y_v$ where $\bar Y_e$ is a compact Hadamard space. 
    
    \item if $\langle t_1\rangle=Z(G_{v_1}), \langle t_2\rangle=Z(G_{v_2})$ then $\langle t_1, t_2\rangle$ is a finite index subgroup of $G_e$.  
\end{enumerate}

We first choose, in a $G$--equivariant way, a plane $F_e \subset Y_e$ for each $e \in T^1$. 

\begin{defn}[Flip CKA action]
\label{defn:flipadaction}
If for each edge $e:=[v,w]\in T^1$,  the   boundary line  $\ell=\bar Y_v \cap F_e$ is parallel  to the $\R$--line in $Y_w = \bar{Y}_w \times \R$, then the CKA action is called \textit{flip}. 
\end{defn}

\subsection{Proof of Theorem~\ref{maintheorem}}
In this section, we are going to prove Theorem~\ref{maintheorem}. We first review some results that will be used.

In \cite{Wis00}, Wise introduces the concept of an \textit{omnipotent group} which has been widely used in subgroup separability. 

\begin{defn}
\label{defn:omnipotence}
A set of group elements $h_1, \cdots, h_r$ in a group $H$ is called {\it{independent}} if whenever $h_i$ and $h_j$ have conjugate powers then $i =j$.
A group $H$ is {\it{omnipotent}} if whenever $\{h_1,\cdots, h_r\}$ ($r\ge 1$) is an independent set of group elements, then there is a positive integer $p\ge 1$ such that for every choice of positive integers $\{n_1,\cdots, n_r\}$, there is a finite quotient $\varphi \colon H\to \hat H$  such that $\varphi(\hat h_i)$ has order $n_ip$ in $\hat H$  for each $i$. 
\end{defn}

It is worth mentioning that free groups \cite{Wis00}, surface groups \cite{Baj07}, Fuchsian groups \cite{Wil10}  and virtually special hyperbolic groups \cite{Wis00} all belong to the omnipotent group category. However, it is a longstanding open question whether every hyperbolic group is residually finite. Wise suggested that if every hyperbolic group is residually finite, then any hyperbolic group would be considered an omnipotent group \cite[Remark~3.4]{Wis00}).

By  \cite[Theorem~II.6.12]{BH99}, each vertex group $G_{\hat{v}}$ of the admissible group $G$ contains a subgroup $K_{\hat{v}}$ intersecting trivially with $Z_{\hat{v}}$ so that the direct product $K_{\hat{v}} \times Z_{\hat{v}}$ is a finite subgroup of $G_{\hat{v}}$. The image of $K_{\hat{v}}$ in the quotient $Q_{\hat{v}} = G_{\hat{v}}/ Z_{\hat{v}}$ is of finite index of $Q_{\hat{v}}$. Since $Q_{\hat{v}}$ is omnipotent and then is residually finite, we can assume that $K_{\hat{v}}$ is torsion-free.

 A collection of finite index subgroups $\{G'_{\hat e}, G'_{\hat v} \, \bigl | \, \hat v \in  V(\Gamma),  \hat e \in E(\Gamma)\}$ of vertex and edge groups of  $G = \pi_1(\mathcal{G})$ is called {\it{compatible}} if $G'_{\bar{\hat e}} = \tau_{e} (G'_{\hat e})$ and whenever $\hat v = \hat e_{-}$ we have
$G'_{\hat e} = G'_{\hat v} \cap  G_{\hat e}$.
When studying the virtual properties of a graph of groups $G$, it is frequently necessary to create a finite index subgroup $G'$ from a set of finite index subgroups of vertex groups. This can be accomplished using the following theorem. 

\begin{thm}\cite[Theorem~7.50]{DK18}
\label{thm:finitecompatible}
Let $G$ be the fundamental group of a  graph of groups $\mathcal{G} = (\Gamma, \{G_{\hat v}\}, \{G_{\hat e} \}, \{\tau_{\hat e} \})$. 
 For every compatible collection $\{G'_{\hat e}, G'_{\hat v} \, \bigl | \, \hat v \in  V(\Gamma),  \hat e \in E(\Gamma)\}$ of $\mathcal{G}$, there exists a finite index subgroup $G' < G$ such that $G' \cap G_{\hat v} = G'_{\hat v}$ and $G' \cap G_{\hat e} = G'_{\hat e}$ for every vertex $\hat v$ and edge $\hat e$.
\end{thm}

\begin{lem} \cite[Lemma~4.8]{HNY23}
   \label{FindexSubgrpsLem}
Let $\{\dot K_{\hat{v}} \le K_{\hat v}: v\in V(\Gamma) \}$ be a   collection of finite index subgroups. Then there exist   finite index subgroups $\ddot K_{\hat v}$ of $\dot K_{\hat v}$, $G'_{\hat e}$ of $G_{\hat e}$ and $Z'_{\hat v}$ of $Z_{\hat{v}}$ so that the  collection of finite index subgroups $\{G'_{\hat e}, G'_{\hat v} =\ddot K_{\hat v} \times Z'_{\hat v}: v\in V(\Gamma), e\in E(\Gamma) \}$ is compatible.
\end{lem}

We are now ready for the proof of Theorem~\ref{maintheorem}. We recall the statement of Theorem~\ref{maintheorem} for the convenience of the reader.

\maintheorem*

\begin{proof}
According to \cite[Lemma~4.6]{NY23}, there is a subgroup of $G$ that has a finite index of at most 2 and is also an admissible group, with a bipartite underlying graph. For simplicity, we still refer to this subgroup as $G$.
Using Lemma~\ref{FindexSubgrpsLem}, we obtain  an admissible group $\mathcal{G'} = \bigl (\Gamma', \{G'_{\hat u} = K_{\hat u} \times \Z_{\hat u} \}, \{G'_{\hat e}  \}, \{\tau_{\hat e} \} \bigr )$ where 
\begin{enumerate}
      \item  $K_{\hat u}$ is a torsion-free, omnipotent CAT(0), nonelementary hyperbolic group. 
    \item $G' : = \pi_1(\mathcal{G}')$ is a finite index subgroup of $G$.
\end{enumerate}
 
 For each vertex $\hat{u}_i \in V(\Gamma')$, 
let $Y_{\hat{u}_i}$ be a CAT(0) hyperbolic space such that $K_{\hat{u}_i} \curvearrowright Y_{\hat{u}_i}$  geometrically. Fix a generator $t_{\hat{u}_i}$ of the factor $\Z_{\hat{u}_i}$ of $K_{\hat{u}_i} \times \mathbb{Z}_{\hat{u}_i}$. Then $G'_{\hat{u}_i} = K_{\hat{u}_i} \times \langle t_i \rangle$   acts geometrically on the $CAT(0)$ space $X_{\hat{u}_i} := Y_{\hat{u}_i} \times \mathbb R$. 
Let $\hat e$ be an oriented edge in $\Gamma'$ such that   $\hat{e}_{-} = \hat{u}_i$. The image $\pi'_{\hat{u}_i}(G'_{\hat e}) \le K_{\hat{u}_i}$ under the projection $\pi'_{\hat{u}_i} \colon G'_{\hat{u}_i} \to K_{\hat{u}_i}$  is an infinite cyclic subgroup generated by an element $k_{\hat e} \in K_{\hat{u}_i}$. The hyperbolic element $k_{\hat e}$ gives rise to a totally geodesic torus $T_{\hat e}$ in the quotient space $X_{\hat{u}_i}/G'_{\hat{u}_i}$ with basis denoted by $([k_{\hat e}], [t_{\hat{u}_i}])$. 
We re-scale $Y_{\hat{u}_i}$ so that the translation length of $k_{\hat e}$ is equal to that of $t_{\hat{u}_i}$ for each $i$. 
Let $$f_{\hat{e}} \colon T_{\hat e} \to T_{\bar{\hat e}}$$ be a \textit{flip} isometry respecting these lengths, that is, an orientation-reversing isometry mapping $[k_{\hat e}]$ to $[t_{\hat{e}_{-}}]$ and $[t_{\hat{e}^{+}}]$ to $[k_{\bar{\hat{e}}}]$.

Let $M$ be the space obtained from taking the disjoint union of compact spaces $\bigsqcup_{\hat{u}_i \in V(\Gamma')} X_{\hat{u}_i}/ G'_{\hat{u}_i}$ and glue these spaces accordingly via isometry $f_{\hat{e}} \colon T_{\hat e} \to T_{\bar{\hat e}}$ with $\hat e$ varies oriented edges on the underlying graph $\Gamma'$.

The fundamental group $\pi_1(M)$ has a graph of groups structure as follows: 
\begin{itemize}
    \item for each vertex $\hat{u}_i$, the associated vertex group is $\pi_1( X_{\hat{u}_i} / G'_{\hat{u}_i})$;
    \item  for each oriented edge $\hat{e}$, the associated edge group is $\pi_1(T_e)$. Edge monomorphisms are $(f_{\hat{e}})_{*} \colon \pi_1(T_{\hat e}) \to \pi_1 (T_{\bar{\hat e}})$ induced by $f_{\hat{e}} \colon T_{\hat e} \to T_{\bar{\hat e}}$.
\end{itemize}

There is a metric on $M$ which makes $M$ into a locally $CAT(0)$ space (see e.g. \cite[Proposition II.11.6]{BH99}).

Let $\widetilde{M} \to M$ be the universal cover of $M$. By the Cartan-Hadamard Theorem, the universal cover $\widetilde M$ with the induced length metric from $M$ is a $CAT(0)$ space, and hence $\pi_1(M)$ is a CAT(0) admissible groups as $\pi_1(M)$ acts geometrically on $\widetilde M$.

As two admissible groups $G'$ and $\pi_1(M)$ have the same underlying graph, same vertex groups, and same edge groups. The only difference is gluing edge maps. We thus can apply  Theorem~\ref{thm:maintheorem1} to conclude that $G'$ and $\pi_1(M)$ are quasi-isometric, and hence  $G$ is quasi-isometric to $\pi_1(M)$ since $G'$ is a finite index subgroup of $G$. 
\end{proof}

Below, we give the proof of Corollary~\ref{cordivergence}. We need several lemmas. \begin{lem} \cite[Corollary 4.17]{BD14}
\label{lem:upperdivergence}
    If a finitely generated group $G$ is strongly thick of order at most $n$, then the divergence of $G$ is bounded above by a polynomial of degree $n$ + 1.
\end{lem}

\begin{lem} \cite[Theorem 6.4]{BD14}
\label{lem:lower-div}
  Let $\gamma$ be a Morse quasi-geodesic in a CAT(0) metric space $X$.
Then the divergence of $X$ is at least quadratic.  
\end{lem}

Suppose that $G$ contains a vertex group of type $\mathcal{H}$. By the normal form theorem, for each connected subgraph $\Gamma'$ of $\Gamma$, there is a subgroup $G_{\Gamma'}\leq G$ which is the fundamental group of the graph of groups with underlying graph $\Gamma'$, and with vertex, edge groups, and edge monomorphisms coming from $\cG$.
Let $\Lambda$ be the full subgraph of $\Gamma$ with vertex set $\{v\in V\Gamma : G_v \text{ is type $\mathcal S$}\}$. For each component $\Gamma'$ of $\Lambda$, we say that $G_{\Gamma'}$ is \begin{enumerate}
	\item a \emph{maximal admissible component} if $\Gamma'$ contains an edge;
	\item an \emph{isolated type $\mathcal S$ vertex group} if  $\Gamma'$ consists  of a single vertex of type $\mathcal S$.
\end{enumerate}  Recall that if $G_v$ is a vertex group of type $\mathcal{H}$, then it is a relatively hyperbolic group to $\PP_v$.

We remark that every graph of groups is obtained by iterating amalgamated products and HNN extensions. By applying the Combination Theorem of relatively hyperbolic groups \cite[Theorem~0.1]{Dah03} to our setting $\mathcal{G}'$, specifically (2) and (3) of \cite[Theorem~0.1]{Dah03} for amalgamated products and (4) for HNN extensions, we obtain the following:

\begin{lem} \cite[Lemma~4.1]{Ngu25}
\label{lem:RHGstructure}
    Let $G_1, . . . , G_k$ be the maximal admissible components
and isolated vertex groups of type $\mathcal{S}$ of an extended admissible group $G$. Let $G_{e_1}, \ldots, G_{e_m}$ be the edge groups so that both its associated vertex groups $G_{(e_i)_{\pm}}$ are of type $\mathcal{H}$, and let $T_1, \ldots, T_{\ell}$ be groups in $\cup \PP_v$ which are not edge groups of $G$.
Then $G$ is hyperbolic relative to 
\[
\PP = \{G_i\}_{i=1}^{k} \cup \{G_{e_s}\}_{s =1}^{m} \cup \{T_i\}_{i=1}^{\ell}
\]
\end{lem}

\cordivergence*

\begin{proof}

We consider the following two cases.

\textbf{Case~1:} $G$ contains no vertex group of type $\mathcal{H}$. In this case $G$ is an admissible group. We first show that the upper bound of the divergence is quadratic.  By \cite[Corollary~3.11]{MN24}, the inclusion of a vertex group $G_v \to G$ is a quasi-isometric embedding, and hence for any two points $x, y \in G_v$, a geodesic $\gamma$ in $G_v$ connecting $x$ to $y$ will be a uniform quasi-geodesic in $G$. This shows that the graph $G_v$ satisfies the quasi-convexity property as defined in~\cite[\S 4.1]{BD14}.
	Since every asymptotic cone of a vertex group of $G$ is without cut-points, it follows that vertex groups of $G$ are strongly algebraically thick of order zero in the sense of~\cite{BD14}. We have that $G$ is strongly thick of order at most $1$ since a graph of groups with infinite edge groups and whose vertex groups is thick of order $n$, is thick of order at most $n+1$, by~\cite[Proposition~4.4 \& Definition 4.14]{BD14}.
Using Lemma~\ref{lem:upperdivergence}, we have that the divergence of $G$ is at most quadratic.

Now we consider the lower bound of the divergence.
According to Theorem~\ref{maintheorem}, there exists a CAT(0) admissible group $G'$ so that $G$ are $G'$ are quasi-isometric. Pick any infinite order group element $g \in G'$ which is not conjugate into any vertex group of $G'$. Then by \cite[Corollary 6.16]{NY23}, $g$ is a Morse element in $G'$. According to Lemma~\ref{lem:lower-div}, the divergence of $G'$ is at least quadratic, and hence the divergence of $G$ must be at least quadratic since   divergence is a quasi-isometric invariant \cite{Ger94}. Therefore the divergence of $G$ is quadratic.

\textbf{Case~2:} $G$ contains at least one vertex group of type $\mathcal{H}$. In this case  $G$ has the natural relatively hyperbolic structure described by Lemma~\ref{lem:RHGstructure}. According to \cite[Theorem~1.3]{Sis12}, the divergence of a relatively hyperbolic group is exponential, and hence the divergence of $G$ is exponential.
\end{proof}

\corQRboundary*

 \begin{proof}
    We consider the following two cases.

    \textit{Case~1:} $G$ contains no vertex group of type $\mathcal{H}$. In this case, we apply Theorem~\ref{maintheorem} to obtain a CAT(0) admissible group $G'$ so that $G$ and $G'$ are quasi-isometric. In \cite[Section 5]{NQ25}, the authors prove that for the CAT(0) admissible group $G'$, its quasi-redirecting boundary is well-defined. Since quasi-redirecting boundary is a quasi-isometric invariant, it follows that $G$ has well-defined quasi-redirecting boundary.

    \textit{Case~2:} $G$ contains at least one vertex group of type $\mathcal{H}$. Let $\PP = \{G_i\}_{i=1}^{k} \cup \{G_{e_s}\}_{s =1}^{m} \cup \{T_i\}_{i=1}^{\ell}$ be the peripheral subgroups of $G$ in the relatively hyperbolic structure of $G$ given by Lemma~\ref{lem:RHGstructure}.
Each $T_i$ and $G_{e_s}$ are quasi-isometric to $\mathbb R^2$ so they do have well-defined quasi-redirecting boundary. The well-defined quasi-redirecting boundary of each $G_i$ is confirmed by Case~1. Thus   we have shown that each peripheral subgroup in this relatively hyperbolic structure has well-defined quasi-redirecting boundary, and hence it follows from \cite[Theorem~D]{NQ25} that $G$ has well-defined quasi-redirecting boundary.
 \end{proof}

\freebycyclic*

 \begin{proof}
   It is shown in the proof of \cite[Lemma~6.8]{BGGH25} that there is a finite index subgroup $\Gamma$ of  $G$ such that $\Gamma$ is an admissible group $\Gamma$ where each vertex group $\Gamma_v$ of $\Gamma$ is a direct product of a free group $F_v$ with $\Z$. Note that $\Gamma_v$ is a CAT(0) group and $F_v$ is omnipotent.  It follows from Corollarry~\ref{corQRboundary} that $\Gamma$ has well-defined quasi-redirecting boundary. Since quasi-redirecting boundary is well-behaved under quasi-isometries, it follows that $G$ has well-defined quasi-redirecting boundary as $\Gamma$ has finite index in $G$.
 \end{proof}

\section{Subgroups of extended admissible groups}
\label{sec:subgroups}
In this section, we study various aspects of subgroups of extended admissible groups.

\begin{prop}
\label{prop:Rapiddecay}Let $G$ be an extended admissible group. Then $G$ has Rapid Decay property. As a consequence, the following groups can not be embedded in extended admissible groups.
\begin{itemize}
    \item Amenable groups with exponential growth.
    \item The Baumslag–Solitar group 
$
B S(p, q)=\left\langle a, b \mid b a^p b^{-1}=a^q\right\rangle
$ where $p\ne\pm q,\,p,q\in \mathbb{Z}$ are nonzero integers.
\item Thompson’s group, $SL_{n}(\mathbb Z)$ with $n \ge 3$, intermediate growth groups.
\end{itemize}
\end{prop} 
\begin{proof} We consider the following two cases.

\textit{Case~1:} $G$ contains no vertex group of type $\mathcal{H}$. In this case, each vertex group is a central extension of a hyperbolic group, hence it is Rapid Decay (\cite{Nos92}). Instead of recalling the precise definition of Rapid Decay ,we refer the reader to \cite{Cha17} for a clear survey on this property, since we only require some of its basic properties.  Also, both vertex groups and edge groups are quasi-isometric embedded \cite[Lemma~2.6, Corollary~3.11]{MN24}, hence these groups are undistorted in $G$ and they have loose polynomial distortion \cite{CG24}. Applying \cite[Proposition~1.3]{CG24}, $G$ has Rapid Decay property. 

\textit{Case~2:} $G$ contains a vertex group of type $\mathcal{H}$. In this case, $G$ is hyperbolic relative to $$\mathbb{P}=\{G_i\}_{i=1}^k\cup \{G_{e_s}\}_{s=1}^m\cup \{T_i\}_{i=1}^l,$$ which is shown in Lemma~\ref{lem:RHGstructure}. Each vertex group $G_v$ of type $\mathcal{H}$ is relative to $\mathbb{P}_v$, and each group in $\mathbb{P}_v$ is Rapid Decay \cite[Theorem~3.1.7]{Jol90}, hence $G_v$ is also Rapid Decay by \cite[Theorem~1.1]{DS05}. Since $G$ is relatively hyperbolic to $\mathbb{P}$, by applying \cite[Theorem~1.1]{DS05} again, $G$ has the Rapid Decay property. 

We remark here that the Rapid Decay property is preserved by passing to subgroups. In \cite{Cha17}, the author lists several classes of groups which do not have the Rapid Decay property, including amenable groups with exponential growth, Baumslag–Solitar groups, Thompson's group, $SL_{n}(\mathbb{Z})$ with $n \ge 3$, and groups of intermediate growth. %Since we have shown that extended admissible groups have the Rapid Decay property, it follows that these mentioned groups cannot be embedded in extended admissible groups.
\end{proof}

\begin{defn}
Let $G$ be a group. A subgroup $H$ of $G$ is {\it separable} if and only if for all $g \in G \backslash H$, there exists a finite index subgroup $K \leq G$ such that $H \leq K \leq G$ and $g \not \in K$. The group $G$ is called {\it locally extended residually finite} (LERF) if any finitely generated subgroup of $G$ is separable.
\end{defn}

%\begin{defn}
%If $\Gamma' \subseteq \Gamma$ is connected, then {\it the graph of groups carried} by $\Gamma'$ is the graph of groups $\mathcal{G}'$ with underlying graph $\Gamma'$ such that the vertex $v \in V(\Gamma') \subseteq V(\Gamma)$ is labelled by $G_v$, and the edge $e \in E(\Gamma') \subset E(\Gamma)$ is labelled $G_e$ (with obvious edge maps). There is a natural map $G' \rightarrow G$, where $G = \pi_1 (\mathcal{G})$ and $G' = \pi_1(\mathcal{G}')$. This map is an injection by the Normal Form Theorem.
%\end{defn}

We consider the following \textit{Croke--Kleiner group}:
\[
L = \langle i, j, k, l \mid [i,j], [j,k], [k,l] \rangle .
\]
The group $L$ is the fundamental group of a torus complex (see Example~\ref{exmpleextended}). 
Note that it is also the fundamental group of a graph manifold and this is a right-angled Artin group on the line graph with four vertices and three edge. 
This group appears in \cite{CK00} as an example of  admissible groups used to show the existence of a group acting geometrically on distinct CAT$(0)$ spaces whose visual boundaries are not homeomorphic (see \cite{CK00}). 
This observation has served as motivation for several boundary constructions in recent years.

In \cite{NW01}, the authors prove the following.

\begin{lem}\cite[Theorem~1.2]{NW01} The Croke-Kleiner group $L$ is not LERF.\label{lem:NW} 
\end{lem}

%Consider the group $G = (F_{m} \times \Z)*_{A} (F_{n} \times \Z)$ where $n, m \ge 2$ and where the amalgamated subgroup $A$ is a maximal $\Z^2$--subgroup in each factor. Suppose that the embedding of the amalgamated subgroup does not identify the centers of the factors. Then there is an embedding $L \to G$.

%\begin{lem}[Gromov] \label{lem:Pingpongcorollary} Let $\Gamma$ be a $\delta$--hyperbolic group. If $u$ and $v$ are two elements in $\Gamma$ such that $u \cdot v \neq v \cdot u$, then for all sufficiently large $m$ and $n$ we have $\langle u^{m}, v^{n} \rangle  \cong F_{2}$ \end{lem}

\begin{prop}
Suppose that an extended admissible group $G$ contains at least one maximal admissible component then there is an embedding $L \to G$. In particular, $G$ is not LERF.
\end{prop}

\begin{proof}
It suffices to consider $G$ as an admissible group, since any group containing a non-LERF subgroup is itself not LERF.

According to \cite[Lemma~4.6]{NY23}, there is a subgroup $\dot G$  of $G$ that has a finite index of at most 2 and is also an admissible group, with a bipartite underlying graph. We denote the graph of groups structure of $\dot G$ by $\mathcal{K} = (\Gamma, \{\dot G_{\hat v}\}, \{\dot G_{\hat e} \}, \{\tau_{\hat e} \})$ where $\tau_{\hat e}$ is an isomorphism $\dot G_{\hat e} \to \dot G_{\bar {\hat e}}$.
% We are going to show that $\dot G$ contains an amalgamated free product $(F_{2} \times \Z)*_{A}(F_{2} \times \Z)$ where $A$ is a maximal $\Z^2$ subgroup in each factor. By Lemma~\ref{lem:NW}, the amalgamated free product $(F_{2} \times \Z)*_{A}(F_{2} \times \Z)$ is not LERF, and thus $G$ is not LERF since a subgroup of a LERF--group is LERF.
%We consider the following cases:

Pick an edge $\hat e$ in the underlying graph of the admissible group $\dot G$ with two distinct vertices $\hat v = \hat e_{-}$ and $\hat w =\hat e_{+}$. 
Choose a generator $\xi_{\hat v}$ of $Z_{\hat v} : = Z(\dot G_{\hat v}) \cong \Z$, and choose a generator $\xi_{\hat w}$ of $Z_{\hat w} : = Z(\dot G_{\hat w}) \cong \Z$. We recall that the subgroup generated by $Z_v$ and $\tau_{\bar{\hat e}} (Z_w)$ has finite index in $G_{e} \cong \Z^2$ (see (4) in Definition~\ref{defn:extended}). Since $\xi_{\hat v}$ is not contained in $Z_{\hat w}$, there exists an element $h_{\hat{w}}$ in $G_{\hat w}$ so that $\xi_{\hat v}$ does not commute with $h_{\hat w}$. Similarly, there exists an element $h_{\hat{v}}$ in $G_{\hat v}$ so that $\xi_{\hat w}$ does not commute with $h_{\hat v}$. Since $\xi_{\hat v} \in Z_{\hat v}$ and $\xi_{\hat w} \in Z_{\hat w}$, we have that $[\xi_{\hat v}, h_{\hat v}] = 1$, $[\xi_{\hat w} , h_{\hat w}] =1$ and $[\xi_{\hat v}, \xi_{\hat w}] = 1$. Let $\check G : = \langle h_{\hat v} , \xi_{\hat v}, \xi_{\hat w}, h_{\hat v} \rangle $ be the subgroup of $\dot G$. We consider the map $\psi \colon L \to \check{G}$ given by 
\[
i \mapsto h_{\hat v}, j \mapsto \xi_{\hat v}, k \mapsto \xi_{\hat w}, l \mapsto h_{\hat w}
\]
Since $\psi(r) =1$ for every relator $r$ in  $L$, it follows that $\psi$ is a homomorphism. Normal forms show that the homomorphism $\psi$ is injective. Thus there is an embedding $L \to G$. Since  $L$ is not LERF (see Lemma~\ref{lem:NW}), it follows that $G$ is not LERF.
\end{proof}

\begin{defn}
\label{defn:malnormalheight}
    Recall that a subgroup $H \le G$ is \textit{malnormal} if $H \cap g H g^{-1}$ is
trivial for all $g \notin H$, and is \textit{almost malnormal} if $H \cap g H g^{-1}$ is finite for all $g \notin H$. Let $H \le G$. The \textit{height} of $H$ in $G$ is the largest number $n \ge 0$ so that there are $n$ distinct cosets $\{ g_1 H, g_2 H, \ldots, g_{n} H \}$ so that the intersection of
conjugates $g_{i} H g_{i}^{-1}$ is infinite. Thus finite groups have height $0$, infinite almost malnormal subgroups have height $1$, and so on.
\end{defn}

\begin{defn}
\label{defn:quasiconvex}
 Let $G$ be a finitely generated group and $H$ a subgroup of $G$.  We say $H$ is \emph{strongly quasiconvex} in $G$ if for any $L \ge 1$, $C \ge 0$ there exists $M = M(L,C)$ such that every $(L,C)$--quasi-geodesic in $G$ with endpoints in $H$ is contained in the $M$--neighborhood of $H$. 
\end{defn}

In \cite{Tra19}, Tran shows that strongly quasiconvex subgroups in a finitely generated group have finite height. While the equivalence of strong quasiconvexity and finite height has been established for extended admissible groups in \cite{Ngu24}, the relationship with virtual malnormality in the context of extended admissible groups has not been explicitly treated.

In the setting of relatively hyperbolic groups, Hruska-Wise in \cite{HW09} prove the following result mentioning that this result is new even in the hyperbolic case. Here we recall a subgroup $H$ of $G$ is \textit{separable} if and only if for all $g \in G \backslash H$, there
exists a finite-index subgroup $K \le G$ such that $H \le K \le G$ and $g \in K$.

\begin{prop} \cite[Theorem~9.3]{HW09}
Let $H$ be a separable, relatively quasiconvex subgroup of the relatively
hyperbolic group $G$. Then there is a finite index subgroup $K$ of $G$ containing $H$ such
that $H$ is relatively malnormal in $K$.      
\end{prop}

We generalize this result to a broader setting by showing that strongly quasiconvex and separable subgroups are virtually almost malnormal. This result applies to the setting of extended admissible groups and may be of independent interest.

\separable*

\begin{proof}
We first claim that there are only finitely many double cosets $Hg_1H,Hg_2H, \cdots, Hg_nH$ such that $H \cap g_i Hg_i^{-1}$ is infinite. Indeed, suppose $\{g_i|i\in I\}$ is a collection of cosets such that $H\cap g_iHg_i^{-1}$ is infinite for each $i$. We fix a finite generating set $S$ of $G$. By the proof of Theorem 4.15 in \cite{Tra19} there is a constant $C$ such that $d_S(H,g_iH)<C$ for each $i$. Thus we can translate $g_iH$ by an element of $H$ to obtain a coset $hg_iH$ intersecting the ball of radius $C$ in the Cayley graph $\Gamma(G,S)$ centered at the identity. Since this ball is finite, it follows that the cosets $g_iH$ lie in only finitely many double cosets $Hg_iH$. 

 Since $H$ is separable, there exists a finite index subgroup $K$ of $G$ containing $H$ and $g_i\notin K$ for each $i$. If $k\in K-H$ and $H \cap kHk^{-1}$ is infinite, then $kH = hg_iH$ for some $g_i$ and some $h \in H$. Also $H$ is a subgroup of $K$. Therefore, $g_i$ is a group element in $K$, contradicting our choice of $K$. Consequently $H$ is almost malnormal in $K$.

 Suppose that $G$ is virtually torsion-free. $G_1 < G$ be the torsion free finite index subgroup of $G$ and let $G_2 = K \cap G_1$. We then have $H_2 = H \cap G_2$ malnormal in $G_2$ using \cite[Lemma~4.24]{AGM16}.

 Now we assume that $G$ is the fundamental group of a graph 3-manifold $M$. 
    If $H$ is virtually malnormal then $H$ is virtually finite height and hence $H$ is strongly quasiconvex in $\pi_1(M)$ by \cite{NTY21}. Now we assume that $H$ is strongly quasiconvex in $\pi_1(M)$. By \cite{NS20} $H$ must be separable in $\pi_1(M)$ since otherwise the distortion of $H$ in $\pi_1(M)$ is quadratic or exponential which contradicts to the fact $H$ is strongly quasiconvex in $\pi_1(M)$.Therefore $H$ is virtually malnormal in $\pi_1(M)$.
\end{proof}

\bibliographystyle{alpha}

%\bibliography{hoang}
\end{document}